\documentclass[a4pitaper,11pt]{amsart}
\usepackage{amsmath,amsthm,amssymb}
\usepackage[mathscr]{eucal}
 \usepackage{cite}
\usepackage{upgreek}
\usepackage[bookmarks,bookmarksnumbered]{hyperref}
\usepackage{enumerate}

\numberwithin{equation}{section}

\newcommand{\car}{\curvearrowright}

\theoremstyle{plain}
\newtheorem{main}{Theorem}

\newtheorem{mprop}[main]{Proposition}

\newtheorem{theorem}{Theorem}[section]

\newtheorem{lemma}[theorem]{Lemma}
\newtheorem{proposition}[theorem]{Proposition}

\theoremstyle{definition}

\newtheorem{definition}[theorem]{Definition}
\newtheorem{example}[theorem]{Example}

\newtheorem{remark}[theorem]{Remark}

\begin{document}

\title[Solid ergodicity and orbit equivalence rigidity for coinduced actions]
{Solid ergodicity and orbit equivalence rigidity for coinduced actions}

\author[D. Drimbe]{Daniel Drimbe}
\address{Department of Mathematics, University of Regina, 3737 Wascana Pkwy, Regina, SK S4S 0A2, Canada.}
\email{daniel.drimbe@uregina.ca}
\thanks{The author was partially supported by PIMS fellowship.}

\begin{abstract} 
We prove that the solid ergodicity property is stable with respect to taking coinduction for a fairly large class of coinduced action. More precisely, assume that $\Sigma<\Gamma$ are countable groups such that $g\Sigma g^{-1}\cap \Sigma$ is finite for any $g\in\Gamma\setminus\Sigma$. Then any measure preserving action $\Sigma\car X_0$ gives rise to a solidly ergodic equivalence relation if and only if the equivalence relation of the associated coinduced action $\Gamma\car X$ is solidly ergodic. We also obtain orbit equivalence rigidity for such actions by
showing that the orbit equivalence relation of a rigid or compact measure preserving action $\Sigma\car X_0$ of a property (T) group is ``remembered'' by the orbit equivalence relation of $\Gamma\car X$.

\end{abstract}

\maketitle

\section{Introduction}
Any countable p.m.p. equivalence relation $\mathcal R$ can be written as the {\it orbit equivalence relation} $\mathcal R_{\Gamma\car X}=\{(x_1,x_2)\in X\times X| \Gamma\cdot x_1=\Gamma \cdot x_2\}$ associated to a p.m.p. action $\Gamma\car (X,\mu)$ of a countable group $\Gamma$ on a standard probability space $(X,\mu)$ \cite{FM75}. 
The study of countable equivalence relations is a central theme in orbit equivalence and measured group theory and many spectacular innovations have been made in the last two decades, see the surveys \cite{Sh04,Fu09,Ga10}. The goal of the present work is to investigate the equivalence relations that are associated to coinduced actions. More precisely, if we let $\Sigma<\Gamma$ be countable groups, then for any p.m.p. action $\Sigma\car (X_0,\mu_0)$ we can construct in a natural way a p.m.p. action $\Gamma\car (X,\mu)$ called the associated {\it coinduced action}, see Definition \ref{def}.
This paper is motivated by the following question: to what extent does $\mathcal R_{\Gamma\car X}$ remember some information about $\mathcal R_{\Sigma\car X_0}$?

Our first main goal is to investigate the {\it solid ergodicity} property for 
coinduced actions. In Theorem \ref{A:solid} we show that solid ergodicity is stable with respect to taking coinduction for large classes of groups $\Sigma<\Gamma$. To put this into context, we recall that
a major breakthrough in the classification problem of p.m.p. actions and their equivalence relations is the use of S. Popa's deformation/rigidity theory \cite{Po07} within the framework of von Neumann algebras. This has led to a remarkable progress in understanding the equivalence relations and von Neumann algebras arising from certain classes of actions, including Bernoulli actions, see the surveys \cite{Va10,Io12a,Io17}.

A notable discovery in this direction was made by I. Chifan and A. Ioana in \cite{CI08} by showing that any orbit equivalence relation $\mathcal R$ associated to a Bernoulli action
is {\it solidly ergodic}\footnote{The terminology is introduced by Gaboriau in \cite[Definition 5.4]{Ga10}.}, meaning: 
for every subequivalence relation $\mathcal S\subset \mathcal R$, there exists a measurable partition $X=\sqcup_{n\ge 0} X_n$ into $\mathcal S$-invariant subsets such that $\mathcal S_{|X_0}$ is amenable and $\mathcal S_{X_n}$ is non-amenable and ergodic for any $n\ge 1$. A key aspect of their proof is Popa's influential spectral gap rigidity principle used within the deformation/rigidity approach.
As noticed in \cite[Section 4.2]{CI08}, Chifan and Ioana's solidity is strongly connected to the work of D. Gaboriau and R. Lyons on von Neumann's problem \cite{GL07}, see also the survey \cite{Ho11}. 

The solid ergodicity property for equivalence relations is the analogue of the solidity property for von Neumann algebras introduced by N. Ozawa \cite{Oz03} (see the paragraph before Theorem \ref{B:solid}).
By using a different method based on
C$^*$-algebraic techniques \cite{Oz03} and the topological amenability property \cite{Oz08}, Ozawa obtained, in particular, that
the orbit equivalence relation $\mathcal R_{SL_2(\mathbb Z)\car\mathbb T^2}$ is solidly ergodic. However, the approach of \cite{CI08} has been successfully followed by R. Boutonnet to show that the orbit equivalence relation of any mixing Gaussian actions is solidly ergodic \cite{Bo12} and by A. Marrakchi to obtain non-amenable solidly ergodic type III equivalence relations \cite{Ma16}.

Before stating our first result, we recall the definition of coinduced actions.


\begin{definition}\label{def}
Let $\Sigma<\Gamma$ be countable groups and denote $I=\Gamma/\Sigma$. Let $\phi:I\to \Gamma$ be a section map and define the associated cocycle $c:\Gamma\times I\to\Sigma$ by the formula $c(g,i)=\phi^{-1}(gi)g\phi(i)$, for all $g\in\Gamma$ and $i\in I.$ \\
For any action $\Sigma\overset{\sigma_0}{\car} (X_0,\mu_0)$ we define the {\it coinduced action} 
$\Gamma\overset{\sigma}\car (X_0,\mu_0)^I$ of $\sigma_0$ by the formula
$$
(\sigma_g(x))_{gi}=(\sigma_0)_{c(g,i)}(x_{i}), \text{ for all } g\in\Gamma,i\in I \text{ and } x=(x_i)_{i\in I}\in X_0^I.
$$ 
\end{definition}

Remark that any coinduced action of a p.m.p. action is again p.m.p. 

\begin{main}\label{A:solid}
Let $\Sigma<\Gamma$ be countable groups such that $ g\Sigma g^{-1}\cap\Sigma$ is amenable for any $g\in\Gamma\setminus \Sigma$. Let $\Sigma{\car} (X_0,\mu_0)$ be any p.m.p. action such that $\mathcal R_{\Sigma\car X_0}$ is solidly ergodic and let $\Gamma\car (X,\mu)$ be the associated coinduced action.

Then $\mathcal R_{\Gamma\car X}$ is solidly ergodic. 
\end{main}

Thus, roughly speaking, Theorem \ref{A:solid} asserts that the orbit equivalence relation $\mathcal R_{\Gamma\car X}$ remembers whether $\mathcal R_{\Sigma\car X_0}$ is solidly ergodic or not. 
We notice that Theorem \ref{A:solid} generalizes the result of Ioana and Chifan \cite{CI08}, since if $\Sigma$ is trivial, then the associated coinduced action is precisely a Bernoulli action.
Note also that our result implies, in particular, that the orbit equivalence relation of any coinduced action from an amenable subgroup is solidly ergodic. We continue by providing more examples for which Theorem \ref{A:solid} applies.

\begin{example}\label{example} The following classes of countable groups $\Sigma<\Gamma$ satisfy that $ g\Sigma g^{-1}\cap\Sigma$ is amenable for any $g\in\Gamma\setminus \Sigma$.

\begin{enumerate}
\item[(i)] Let $\Gamma=\Delta*_{\Sigma_0}\Sigma$ be a non-trivial amalgamated free product, where $\Sigma_0$ is an amenable group.
\item[(ii)] Let $\Gamma=\Delta\wr_I\Sigma$ be a generalized wreath product such that $\Sigma\car I$ has amenable stabilizers.

\item[(iii)] Let $d:\Lambda\to\Lambda\times\Lambda$ be the diagonal embedding of a group $\Lambda$ for which
the centralizer $C_\Lambda(g)$ is amenable for any $g\in\Lambda\setminus\{e\}$. 
For instance, let $\Lambda$ be any bi-exact torsion-free group \cite{Oz03}, or any countable subgroup of $SL_2(\mathbb R)$ (see the first paragraph of the proof of \cite[Corollary B]{DHI16}). 
We take $\Sigma:=d(\Lambda)$ and $\Gamma:=\Lambda\times \Lambda$.

\item[(iv)] Let $\Gamma$ be a countable group and let $\Sigma$ be a hyperbolically embedded subgroup of $\Gamma$ in the sense of \cite{DGO11} (e.g. $\Gamma$ is a non-elementary group that is hyperbolic relative to $\Sigma$). Then $g\Sigma g^{-1}\cap \Sigma$ is finite for any $g\in\Gamma\setminus\Sigma$ by \cite[Proposition 2.10]{DGO11}. Other examples of groups that admit hyperbolically embedded subgroups are Out$(\mathbb F_n)$ and most mapping class groups. We refer the reader to \cite{DGO11} for more details and more examples.

\end{enumerate}

\end{example}

Remark that any non-amenable ergodic subequivalence relation of a solidly ergodic equivalence relation is {\it strongly ergodic} by \cite[Proposition 6]{CI08}, see Section \ref{terminology} for the definition.
In the setting of Theorem \ref{A:solid}, although $\mathcal R_{\Gamma\car X}$ is not solidly ergodic in general, we show in Proposition \ref{strong ergodicity} that it is strongly ergodic.
It is easy to see that strong ergodicity  is a weaker property than solid ergodicity. Indeed, let $\mathcal R$ be a countable ergodic p.m.p. equivalence relation on $(X,\mu)$. By a result of Popa (see \cite[Proposition 7]{Oz03}), if $\mathcal R$ is not strongly ergodic, then there exists a diffuse von Neumann subalgebra $Q\subset L^\infty(X)$ such that $Q'\cap L(\mathcal R)$ is non-amenable, proving that $\mathcal R$ is not solidly ergodic by Proposition \ref{equivalence}. 

Our next result shows that many subequivalence relations of $\mathcal R_{\Gamma\car X}$ are strongly ergodic, where $\Gamma\car (X,\mu)$ is a coinduced action.
First, we recall the notion of relative amenability for groups. Following \cite[Section 2.3]{BV12}, for some $\Delta,\Sigma<\Gamma$ countable groups, we say that $\Delta$ is {\it amenable relative to } $\Sigma$ {\it inside} $\Gamma$ if the left translation action $\Delta$ on $\Gamma/\Sigma$  admits an invariant mean. This is equivalently for $L(\Delta)$ to be amenable relative to $L(\Sigma)$ inside $L(\Gamma)$ in the sense of Ozawa and Popa \cite{OP07} (see Section \ref{section amenable}). 

\begin{mprop}\label{strong ergodicity}
Let $\Sigma<\Gamma$ be countable groups. Let $\Sigma{\car} (X_0,\mu_0)$ be any p.m.p. action and let $\Gamma\car (X,\mu)$ be the associated coinduced action. Let $\Delta<\Gamma$ be a subgroup that is not amenable relative to $\Sigma$ inside $\Gamma$. 

Then the restriction action $\Delta\car (X,\mu)$ has spectral gap. In particular, $\mathcal R_{\Delta\car X}$ is strongly ergodic.

\end{mprop}

Hence, if $\Sigma$ is an amenable group, then Proposition \ref{strong ergodicity} implies that $\Delta\car (X,\mu)$ has spectral gap for any non-amenable subgroup $\Delta<\Gamma$. See Section \ref{terminology} for the definition of spectral gap and Proposition \ref{vNa version} for a variant of Proposition \ref{strong ergodicity} that applies to von Neumann algebras.

We continue by showing a von Neumann algebraic version of Theorem \ref{A:solid} using Ozawa's solidity notion introduced in \cite{Oz03}. Ozawa proved that any von Neumann algebra $L(\Gamma)$ arising from a bi-exact group $\Gamma$ is {\it solid}, meaning: for every subalgebra $Q\subset L(\Gamma)$, there exists a sequence of projections $(q_n)_{n\ge 0}\subset \mathcal Z(Q)$ with $\sum_n q_n=1$ such that $Qq_0$ is amenable and $Qq_n$ is a non-amenable factor for any $n\ge 1$ \cite{Oz03}. Subsequently, other solid von Neumann algebras have been discovered. J. Peterson showed in \cite{Pe06} that if a countable group $\Gamma$ has a proper cocycle into a multiple of the left regular representation, then $L(\Gamma)$ is solid. 
I. Chifan and A. Ioana obtained in \cite{CI08} that $L(\Delta\wr\Gamma)$ is solid whenever $L(\Gamma)$ is a solid II$_1$ factor and $\Delta$ is amenable. I. Chifan and C. Houdayer then proved a similar stability result for free products by showing that $L(\Gamma_1*\Gamma_2)$ is solid if $L(\Gamma_1)$ and $L(\Gamma_2)$ are solid \cite{CH08}.


In our next result, we are essentially showing that the coinduction procedure preserves Ozawa's solidity for the associated von Neumann algebras.

\begin{main}\label{B:solid}
Let $\Sigma<\Gamma$ be countable groups such that $ g\Sigma g^{-1}\cap\Sigma$ is finite for any $g\in\Gamma\setminus \Sigma$. Let $\Sigma\car (X_0,\mu_0)$ be a p.m.p. action and let $\Gamma\car (X,\mu)$ be the associated coinduced action. 

Assume that $L^\infty(X_0)\rtimes\Sigma$ and $L(\Gamma)$ are  solid von Neumann algebras.

Then $L^\infty(X)\rtimes\Gamma$ is solid.
\end{main}

\begin{remark}
If $\Sigma\car (X_0,\mu_0)$ is a p.m.p. action such that $L^\infty(X_0)\rtimes\Sigma$ is solid, then $\mathcal R_{\Sigma\car X_0}$ is solidly ergodic. By using \cite[Theorem 3.2]{CS11}, note that the converse holds if we assume that $\Sigma$ is a bi-exact group (see Lemma \ref{soliderg}). In combination with the results that are mentioned above, Theorem \ref{B:solid} provides new examples of solid von Neumann algebras.

\end{remark}


Next, we return to the framework of orbit equivalence relations and discuss about the classification of p.m.p. actions up to orbit equivalence. 
We recall that two free ergodic p.m.p. actions $\Gamma\car (X,\mu)$ and $\Lambda\car (Y,\nu)$ are {\it orbit equivalent} if 
their associated orbit equivalence relations are isomorphic. 
The result of \cite{OW80} (see also \cite{Dy58,CFW81}) shows that the class of amenable groups manifests the following orbit equivalence {\it flexibility} property: any two free ergodic p.m.p. actions of countable infinite amenable groups are orbit equivalent. 

In sharp contrast, R. Zimmer discovered in \cite{Zi84} that lattices in higher rank simple Lie groups have orbit equivalence {\it rigidity} properties in the following sense: if two free ergodic p.m.p. actions $SL_n(\mathbb Z)\car (X,\mu)$ and $SL_m(\mathbb Z)\car (Y,\nu)$ are orbit equivalent, for $m,n\ge 3$, then $m=n$. S. Popa then showed that the more general class of property (T) groups also manifests rigidity properties in orbit equivalence. More precisely, Popa proved that any two Bernoulli actions of icc property (T) groups that are orbit equivalent (or more generally, the actions have isomorphic von Neumann algebras) must be conjugate \cite{Po03,Po04}. The most extreme case of rigidity in orbit equivalence is when an action $\Gamma\car (X,\mu)$ is {\it orbit equivalent superrigid}, which roughly means that its associated orbit equivalence relation completely remembers the group $\Gamma$ and the action $\Gamma\car (X,\mu)$. A. Furman showed that ``most'' ergodic p.m.p. actions of higher rank lattices, including $SL_n(\mathbb Z)\car \mathbb T^n$ for $n\ge 3$, are orbit equivalent superrigid \cite{Fu99} and S. Popa proved then that any Bernoulli action of a non-amenable property (T) or product group is orbit equivalent superrigid \cite{Po05,Po06}.
Subsequently, several large classes of orbit equivalent superrigid actions were found, see the introduction of \cite{DIP19}.

Our next result provides orbit equivalence rigidity for a large class of coinduced actions of free product groups.

\begin{main}\label{A:converse}
Let $\Sigma_1\car (X_1,\mu_1)$ and $\Sigma_2\car (X_2,\mu_2)$ be some free ergodic p.m.p. actions of countable non-amenable groups.
Let $\Delta$ be any countable group and denote $\Gamma_1=\Sigma_1*\Delta$ and $\Gamma_2=\Sigma_2*\Delta.$ Assume that for any $1\leq i\leq 2$ one of the following holds:
\begin{itemize}

\item $\Sigma_i\car (X_i,\mu_i)$ is a rigid or compact action; $\Sigma_i$ has property (T) or it is measure equivalent to a product of two non-amenable groups.

\item $\Sigma_i\car (X_i,\mu_i)$ is conjugate to a product action of two infinite groups.
\end{itemize}

Then $\Sigma_1\car (X_1,\mu_1)$ and $\Sigma_2\car (X_2,\mu_2)$ are orbit equivalent if and only if their associated coinduced actions $\Gamma_1\car (X_1,\mu_1)^{\Gamma_1/\Sigma_1}$ and $\Gamma_2\car (X_2,\mu_2)^{\Gamma_2/\Sigma_2}$ are orbit equivalent.
\end{main}

The ``only if'' part follows from a well-known result of L. Bowen \cite{Bo10}, see the last paragraph of the introduction for more details. 
See also Theorem \ref{general} for a more general statement that applies, for instance, to coinduced actions of wreath product groups.

In our result, we used Gromov's notion of measure equivalence \cite{Gr91} (see Definition \ref{def:ME}) and Popa's notion of rigidity for actions \cite{Po01}. Recall that 
a p.m.p. action $\Gamma\car (X,\mu)$ is {\it rigid} if the inclusion $L^\infty(X)\subset M:=L^\infty(X)\rtimes\Gamma$ has the {\it relative property (T)}, i.e. if any sequence of unital tracial completely positive maps $\Phi_n: M\to M$ converging to the identity pointwise in $\|\cdot\|_2$, must converge uniformly on the unit ball of $L^\infty(X)$ \cite[Section 4]{Po01}. See also \cite{Io09} for an ergodic theoretic characterization of rigidity.

\begin{example}\label{example rigid}
Any action of the form $\Gamma\car (G/\Lambda,m)$ is rigid, whenever $G$ is a simple Lie group, $\Gamma<G$ is any Zariski dense countable subgroup, $\Lambda<G$ is a lattice and $m$ is the unique $G$-invariant probability measure on $G/\Lambda$ \cite{IS10}.
In particular, our result applies if we consider $n\ge2$, $G=SL_n(\mathbb R)$, $\Lambda=SL_n(\mathbb Z)$ and $\Gamma=SL_n(R)$,
where $R=\mathbb Z[\sqrt d]$, for a square free integer $d\ge 2$ or $R=\mathbb Z[S^{-1}]$, for a finite non-empty set of primes $S$.
It also applies if we consider $n\ge 3$, $G=SL_n(\mathbb R)$, and $\Gamma,\Lambda<G$ any two lattices (e.g. $\Gamma=\Lambda=SL_n(\mathbb Z)$).
  
\end{example}

See \cite[Theorem D]{IS10} for a more general statement regarding Example \ref{example rigid}. 
In addition to these concrete classes of rigid actions, note that any free product group $\Gamma=\Gamma_1*\Gamma_2$, with $|\Gamma_1|\ge 2$ and $|\Gamma_2|\ge 3$, admits uncountably many non-orbit equivalent free ergodic p.m.p. rigid actions \cite{Ga08}.

\begin{remark}
The conclusion of Theorem \ref{A:converse} can be strengthened if we consider profinite (or more generally, compact) actions by using orbit equivalence superrigidity theorems for profinite actions. More precisely,
if $\Sigma_1\car (X_1,\mu_1)$ is an ergodic profinite action of a property (T) group or $\Sigma_1=SL_2(\mathbb Z[S^{-1}])$, for a finite non-empty set of primes $S$, then $\Sigma_1\car (X_1,\mu_1)$ is ``virtually'' conjugate to $\Sigma_2\car (X_2,\mu_2)$ by using \cite{Io08b} and \cite{DIP19}, respectively. The same result also holds for separately ergodic profinite actions of product groups, see \cite{GITD16} for a precise statement.
\end{remark}

To put Theorem \ref{A:converse} into a better perspective, note that
any free ergodic action of a non-amenable free group $\mathbb F_n$ is orbit equivalent to actions of uncountably many non-isomorphic groups \cite[Theorem 2.27]{MS02}. In particular, the actions of free groups are ``far away'' from being orbit equivalent superrigid. However, our result is essentially showing that if $\Sigma*\Delta\car X$ is the coinduced action of a rigid action $\Sigma\car X_0$ of a property (T) group, then the orbit equivalence relation $\mathcal R_{\Sigma\car X_0}$ can be recovered from the orbit equivalence relation $\mathcal R_{\Sigma*\Delta\car X}$. We can also contrast Theorem \ref{A:converse} with the main result of \cite{Dr15} which asserts that the coinduced actions of property (T) and product groups are actually orbit equivalent superrigid.

Moreover, Theorem \ref{A:converse} provides a converse of a result of L. Bowen \cite[Theorem 1.3]{Bo10}. To put this into context, recall that S. Popa proved in \cite{Po05} orbit equivalence superrigidity for Bernoulli actions of property (T) groups and the question whether any two orbit equivalent Bernoulli actions over a non-amenable free group are necessarily conjugate was left open.
Bowen answered this question by showing the following orbit equivalence {\it flexibility} result: any two Bernoulli actions of a free group are orbit equivalent \cite{Bo10}. Actually, he obtained the following more general result.
If $\Sigma_1\car X_1$ and $\Sigma_2\car X_2$ are some free p.m.p. orbit equivalent actions and $\Delta$ is any countable group, then the associated coinduced actions $\Sigma_1*\Delta\car X_1^{(\Sigma_1*\Delta)/\Sigma_1}$ and $\Sigma_2*\Delta\car X_2^{(\Sigma_2*\Delta)/\Sigma_2}$ are orbit equivalent.
 

{\bf Organization of the paper.}
Besides the introduction this paper has five other sections. In Section 2, we recall some preliminaries and prove a few useful lemmas needed in the remainder of the paper. In Section 3, we prove some rigidity results for von Neumann algebras associated to coinduced actions by building upon  \cite{Po03,Io06a,IPV10}. Sections 4 and 5 are devoted to the proofs of our main results.

{\bf Acknowledgements.} I am very grateful to Adrian Ioana for asking if one can prove a converse to \cite[Theorem 1.3]{Bo10} and for his helpful suggestions. I would also like to thank the anonymous referees for their many comments that greatly improved the exposition of the paper. In particular, I would like to thank them for pointing out an error in a previous draft and for showing me that Proposition B holds under the current generality.

\section{Preliminaries}

\subsection{Terminology}\label{terminology}
In this paper we consider {\it tracial von Neumann algebras} $(M,\tau)$, i.e. von Neumann algebras $M$ equipped with a faithful normal tracial state $\tau: M\to\mathbb C.$ This induces a norm on $M$ by the formula $\|x\|_2=\tau(x^*x)^{1/2},$ for all $x\in M$. We will always assume that $M$ is a {\it separable} von Neumann algebra, i.e. the $\|\cdot\|_2$-completion of $M$ denoted by $L^2(M)$ is separable as a Hilbert space.
We denote by $\mathcal Z(M)$ the {\it center} of $M$, by $\mathcal U(M)$ the {\it unitary group} of $M$ and by $(M)_1$ its unit ball. All inclusions $P\subset M$ of von Neumann algebras are assumed unital. We denote by $e_P: L^2(M)\to L^2(P)$ the orthogonal projection onto $L^2(P)$ and by $E_{P}:M\to P$ the unique $\tau$-preserving {\it conditional expectation} from $M$ onto $P$.

For a von Neumann subalgebra $P\subset M$, we denote by $P'\cap M=\{x\in M|xy=yx, \text{ for all } y\in P\}$ the {\it relative commutant} of $P$ in $M$ and by $\mathcal N_{M}(P)=\{u\in\mathcal U(M)|uPu^*=P\}$ the {\it normalizer} of $P$ in $M$. The {\it quasi-normalizer} of $P$ inside $M$ is the weak closure of $\{x\in M| \text{ there exist } x_1, \dots, x_n,$ $y_1,\dots,y_n$ such that $xP\subset\sum_{i=1}^n Px_i$ and $Px\subset\sum_{i=1}^n y_iP\}.$
We say that $P$ is {\it regular} in $M$ if the von Neumann algebra generated by $\mathcal N_M(P)$ equals $M$ and $P\subset M$ is {\it Cartan} if it is a maximal abelian regular subalgebra.
For two von Neumann subalgebras $P,Q\subset M$, we denote by $P\vee Q$ the von Neumann algebra generated by $P$ and $Q$. {\it Jones' basic construction} of the inclusion $P\subset M$ is defined as the von Neumann subalgebra of $\mathbb B(L^2(M))$ generated by $M$ and $e_P$, and is denoted by $\langle M,e_P \rangle$.

Let $\mathcal R$ be a countable p.m.p. equivalence relation on $(X,\mu)$. We say that $\mathcal R$ is {\it ergodic} if for every measurable set $A\subset X$ satisfying $\mu(\alpha(A)\Delta A)=0$ for all $\alpha\in [\mathcal R]$, we have $\mu(A)\in\{0,1\}$. We say that $\mathcal R$ is {\it strongly ergodic} if for every sequence of measurable sets $A_n\subset X$ satisfying lim$_{n\to \infty} \mu(\alpha(A_n)\Delta A_n)=0$ for all $\alpha\in [\mathcal R]$, we have lim$_{n\to\infty} \mu(A_n)(1-\mu (A_n))=0$.
Note that a p.m.p. action $\Gamma\car (X,\mu)$ is ergodic (or, strongly ergodic) if and only if $\mathcal R_{\Gamma\car X}$ is.
We say that $\Gamma\car (X,\mu)$ has spectral gap if the Koopman representation of $\Gamma$ on $L^2(X,\mu)\ominus \mathbb C1$ has no invariant vectors and notice that spectral gap implies strong ergodicity.

Finally, if $\Gamma\car I$ is an action on a non-empty set $I$ and $F\subset I$ is a subset, we denote by Stab$(F)=\{g\in\Gamma|gi=i, \text{ for all }i\in F\}$ the stabilizer and by Norm$(F):=\{g\in\Gamma| gF=F\}$ the normalizer of $F$. Note that if $F$ is finite, then Stab$(F)$ is a finite index subgroup of Norm$(F)$.

\subsection {Intertwining-by-bimodules} We next recall from  \cite [Theorem 2.1 and Corollary 2.3]{Po03} the powerful {\it intertwining-by-bimodules} technique of S. Popa.

\begin {theorem}[\!\!\cite{Po03}]\label{corner} Let $(M,\tau)$ be a tracial von Neumann algebra and $P\subset pMp, Q\subset qMq$ be von Neumann subalgebras. Let $\mathcal G\subset\mathcal U(P)$ be a subgroup such that $\mathcal G''=P$.

Then the following are equivalent:

\begin{enumerate}

\item There exist projections $p_0\in P, q_0\in Q$, a $*$-homomorphism $\theta:p_0Pp_0\rightarrow q_0Qq_0$  and a non-zero partial isometry $v\in q_0Mp_0$ such that $\theta(x)v=vx$, for all $x\in p_0Pp_0$.

\item There is no sequence $(u_n)_n\subset\mathcal G$ satisfying $\|E_Q(xu_ny)\|_2\rightarrow 0$, for all $x,y\in M$.

\end{enumerate}
\end{theorem}

If one of these conditions holds true, then we write $P\prec_M Q$, and say that {\it a corner of} $P$ {\it embeds into} $Q$ {\it inside} $M$. If $Pp'\prec_M Q$ for every non-zero projection $p'\in P'\cap pMp$, then we write $P\prec_M^s Q.$

We continue with some elementary lemmas. Recall first that a von Neumann algebra $M$ is called {\it properly non-amenable} if $pMp$ is non-amenable for any non-zero projection $p\in M$ (see also the next subsection).
 
\begin{lemma}\label{solid}
Let $Q\subset qMq$ be tracial von Neumann algebras such that $Q$ is properly non-amenable and has diffuse center.

If $N\subset M$ is a von Neumann subalgebra such that $Q\prec_M N$, then $N$ is not solid.

\end{lemma}

{\it Proof.} 
The assumption implies that there exist projections $q_1 \in Q$, $p\in N$, a non-zero partial isometry $v\in pMq_1$ and a one-to-one $*$-homomorphism such that $\theta(x)v=vx$ for any $x\in q_1Qq_1$. Since $Q$ is properly non-amenable, it follows that $P:=\theta (q_1Qq_1)\subset pNp$ is non-amenable and it has diffuse center $R:=\theta (\mathcal Z(Q)q_1)$. Remark that $\tilde R:=R\oplus pNp$ is a diffuse subalgebra of $N$ for which its relative commutant $\tilde R'\cap N$ is non-amenable. Hence, $N$ is not solid. 
\hfill$\blacksquare$

\begin{lemma}\label{small}
Let $N\subset M$ be tracial von Neumann algebras. Let $P\subset pNp$ and $Q\subset qNq$ be von Neumann subalgebras such that $Q\subset qMq$ is regular.
If $P\prec_M Q$, then $P\prec_N Q.$
\end{lemma}

{\it Proof.} If the contrary holds, then there exists a sequence of unitaries $(u_n)_n\subset \mathcal U(P)$ such that $\|E_{Q}(xu_ny)\|_2\to 0$, for any $x,y\in N.$ Thus, $\|E_{Q}(u_ny)\|_2=\|E_{Q}(u_nE_{N}(y))\|_2\to 0$, for any $y\in M.$ Since $Q$ is regular in $qMq$, we obtain that $\|E_{Q}(xu_ny)\|_2\to 0$, for any $x,y\in M$, contradiction.
\hfill$\blacksquare$


\subsection{Relative amenability.}\label{section amenable} 
A tracial von Neumann algebra $(M,\tau)$ is {\it amenable} if there exists a positive linear functional $\Phi:\mathbb B(L^2(M))\to\mathbb C$ such that $\Phi_{|M}=\tau$ and $\Phi$ is $M$-{\it central}, meaning $\Phi(xT)=\Phi(Tx),$ for all $x\in M$ and $T\in \mathbb B(L^2(M))$. By Connes' striking classification of amenable factors \cite{Co76}, $M$ is amenable if and only if $M$ is approximately finite dimensional.
 
We continue by recalling the notion of relative amenability introduced by 
Ozawa and Popa in \cite{OP07}. Let $(M,\tau)$ be a tracial von Neumann algebra. Let $p\in M$ be a projection and $P\subset pMp,Q\subset M$ be von Neumann subalgebras. Following \cite[Definition 2.2]{OP07}, we say that $P$ is {\it amenable relative to $Q$ inside $M$} if there exists a positive linear functional $\Phi:p\langle M,e_Q\rangle p\to\mathbb C$ such that $\Phi_{|pMp}=\tau$ and $\Phi$ is $P$-central. 
Note that $P$ is amenable relative to $\mathbb C$ inside $M$ if and only if $P$ is amenable.

We say that $P$ is {\it strongly non-amenable relative to} $Q$ if $Pp'$ is non-amenable relative to $Q$ for any non-zero projection $p'\in P'\cap pMp$. Note that if $P\subset pMp$ is strongly non-amenable (relative to $\mathbb C$), then $P$ is {properly non-amenable} by \cite[Lemma 2.6]{DHI16}. 


\begin{proposition}\label{splits}
Let $\Gamma\car A$ and $\Gamma\car B$ be some trace preserving actions of a countable group $\Gamma$.
Let $Q\subset qAq$ be a von Neumann subalgebra such that $Q'\cap q(A\rtimes\Gamma)q$ is amenable relative to $A$.

Then $Q'\cap q((A\bar\otimes B)\rtimes\Gamma)q$ is amenable relative to $A\bar\otimes B$.
\end{proposition}

{\it Proof.} The proof is inspired by the proof of \cite[Proposition 3.2]{PV12}. Denote $\mathcal M=(A\bar\otimes B)\rtimes\Gamma$, $M=A\rtimes\Gamma$ and $N=B\rtimes\Gamma$.
Define the $*$-homomorphism $\Delta: \mathcal M\to M\bar\otimes N$, by letting $\Delta((a\otimes b) u_g)=au_g\otimes bu_g$, for all $a\in A, b\in B$ and $g\in\Gamma$. Remark that the assumption implies that $\Delta(Q'\cap q\mathcal Mq)\subset (Q'\cap qMq)\bar\otimes N$ is amenable relative to $A\otimes N$ inside $M\bar\otimes N.$ Therefore, there exists a $\Delta(Q'\cap q\mathcal M q)$-central positive linear functional $\Phi: \Delta(q)\langle M\bar\otimes N, e_{A\otimes N} \rangle \Delta(q)\to\mathbb C$ such that the restriction of $\Phi$ to $\Delta(Q'\cap q\mathcal Mq)$ equals the trace on $\Delta(q)(M\bar\otimes N)\Delta (q).$

Since $E_{A\otimes N}\circ \Delta (x)=\Delta\circ E_{A\bar\otimes B}(x)$ for every $x\in\mathcal M$, notice that we can define a one-to-one $*$-homomorphism $\Delta_1:\langle \mathcal M,e_{A\bar\otimes B} \rangle\to \langle M\bar\otimes N, e_{A\otimes N} \rangle$ by letting $\Delta_1(e_{A\bar\otimes B})=e_{A\otimes N}$ and $\Delta_1(x)=\Delta(x)$, for any $x\in \mathcal M.$ Define the positive linear functional $\Psi: q\langle \mathcal M, e_{A\bar\otimes B} \rangle q \to \mathbb C$ by $\Psi(x)=\Phi(\Delta_1(x))$, for every $x\in q\langle \mathcal M, e_{A\bar\otimes B} \rangle q$. Remark that $\Psi$ is $Q'\cap q\mathcal M q$-central and the restriction of $\Psi$ to $q\mathcal M q$ equals the trace on $q\mathcal M q.$
This shows that $Q'\cap q\mathcal Mq$ is amenable relative to $A\bar\otimes B$ inside $\mathcal M$.
\hfill$\blacksquare$

Although later on we will only use the particular case of Proposition \ref{splits} where $A$ and $B$ are amenable, we note that the proof  handles this general statement.

As a consequence, we obtain the following lemma
that will be needed in the proof of Theorem \ref{B:solid}. 

\begin{lemma}\label{soliderg}
Let $\Gamma\car (X\times Y,\mu\times\nu)$ be a diagonal product p.m.p. action and denote $M=L^\infty(X\times Y)\rtimes\Gamma$.
Let $Q\subset qMq$ be a properly non-amenable von Neumann subalgebra such that there exists a diffuse subalgebra $Q_0\subset \mathcal Z(Q)$ satisfying $Q_0'\cap qMq=Q$ and $Q_0\prec_M L^\infty(X)$.

Then $\mathcal R_{\Gamma\car X}$ is not solidly ergodic.

\end{lemma}

{\it Proof.} 
Since $Q_0\prec_M L^\infty(X)$, there exist projections $q_0\in Q_0,r\in L^\infty(X)$, a non-zero partial isometry $w\in rMq_0$ and a $*$-isomorphism $\theta: Q_0q_0\to R\subset L^\infty(X)r$ such that $\theta(x)w=wx$, for any $x\in Q_0q_0.$ Note that Ad$(w)$ is an isomorphism from $w^*wMw^*w$ to $ww^*Mww^*$ that sends $w^*wq_0(Q_0'\cap qMq)q_0w^*w$ onto $ww^*(R'\cap rMr)ww^*$. Since $Q_0'\cap qMq$ is properly non-amenable, it follows that $R'\cap rMr$ is non-amenable. By applying Proposition \ref{splits}, we get that $R'\cap r(L^\infty(X)\rtimes\Gamma)r$ is non-amenable, hence $\mathcal R_{\Gamma\car X}$ is not solidly ergodic.
\hfill$\blacksquare$






\subsection{Coinduced actions}

Coinduced actions for trace preserving actions are defined as in Definition \ref{def} for p.m.p. actions.

\begin{definition}\label{def:coinduce}

Let $\Sigma<\Gamma$ be countable groups and denote $I=\Gamma/\Sigma$. Let $\phi:I\to \Gamma$ be a section map and define the associated cocycle $c:\Gamma\times I\to\Sigma$ by the formula $c(g,i)=\phi^{-1}(gi)g\phi(i)$, for all $g\in\Gamma$ and $i\in I.$ 

Let $\Sigma\overset{\sigma_0}{\car} (A_0,\tau_0)$ be a trace preserving action and denote $(A,\tau)=(A_0,\tau_0)^I$. For each $i\in I$, we denote by $\pi_i: A_0\to A$ the embedding of $A_0$ as the $i$'th tensor factor. We define the {\it coinduced action} $\Gamma\overset{\sigma}\car (A,\tau)$ by the formula
$$
\sigma_g(\pi_i(a))=(\sigma_0)_{c(g,i)}(\pi_{gi}(a)), \text{ for all } g\in\Gamma,i\in I \text{ and } a\in A_0.
$$
\end{definition}

Throughout this section, we assume the setting of Definition \ref{def:coinduce} and denote $M=A\rtimes\Gamma$.

\begin{lemma}\label{inside}

Let $B\subset A$ be a von Neumann subalgebra such that $B\prec_M A_0^F$ for some subset $F\subset I.$

Then there exists $g\in\Gamma$ such that $B\prec_A A_0^{gF}$.

\end{lemma}

{\it Proof.} Assume that $B\nprec_A A_0^{gF}$ for any $g\in\Gamma$.
By applying an idea of the proof of \cite[Theorem 4.3]{IPP05}, we obtain that there exists a sequence of unitaries $(u_n)_n\subset \mathcal U(B)$ such that 
$$
\|E_{A_0^{gF}}(xu_ny)\|_2\to 0, \text{for all }x,y\in A \text{ and }g\in\Gamma. 
$$
We want to show that $\|E_{A_0^F}(xu_ny)\|_2\to 0$ for all $x,y\in M$. By Kaplansky’s density theorem it is enough to assume that $x=u_{g^{-1}}a$ and $y=a'u_h$ for some $a,a'\in A$ and $g,h\in\Gamma$. Thus, $\|E_{A_0^F}(xu_ny)\|_2=\delta_{g,h}\|E_{A_0^{gF}}(au_na')\|_2\to 0$. This shows that $B\nprec_M A_0^F $, contradiction.
\hfill$\blacksquare$

We continue with the following useful lemma extracted from \cite[Theorem 4.2]{IPV10}.

\begin{lemma}\label{normalizer}
Let $F\subset I$ be a subset and let $Q\subset qMq$ be a diffuse von Neumann subalgebra such that $Q\prec_M A_0^F$ and $Q\nprec_M A_0^{G}$, for any proper subset $G\subset F$.  Denote by $P$ the quasi-normalizer of $Q$ inside $qMq$.

Then $P\prec_M A\rtimes$Norm$(F)$.

\end{lemma}

{\it Proof.} Since $Q\prec_M A_0^F$, there exist some projections $q_0\in Q$, $p\in A_0^F$, a non-zero partial isometry $v\in pMq_0$ and a one-to-one $*$-homomorphism $\phi: q_0Qq_0\to pA_0^Fp$ such that $\phi(x)v=vx$, for any $x\in q_0Qq_0.$ Moreover, we can assume that the support of $E_{A_0^F}(vv^*)$ equals $p.$ The assumption on $F$ implies 
\begin{equation}\label{int}
\phi(q_0Qq_0)\nprec_{A_0^F} A_0^G, \text{ for any proper subset } G\subset F.
\end{equation}
Indeed, by assuming the contrary, there exist a proper subset $G\subset F$, some projections $q_1\in q_0Qq_0, r\in A_0^G$, a non-zero partial isometry $w\in rA_0^F\phi(q_1)$ and a $*$-homomorphism $\psi: \phi(q_1Qq_1)\to rA_0^Gr$ such that $\psi(x)w=wx$, for any $x\in \phi(q_1Qq_1)$. Therefore, $\psi\circ\phi: q_1Qq_1\to rA_0^Gr$ satisfies $(\psi\circ\phi) (x)wv=wvx$. Note that $wv\neq 0$, since $wE_{A_0^F}(vv^*)\neq 0$. By replacing $wv$ by the partial isometry from its polar decomposition, we deduce that $Q\prec_M A_0^G$, contradiction.

Note that \eqref{int} allows us to apply \cite[Lemma 4.1.2]{IPV10} and derive that $vP v^*\subset A\rtimes$Norm$(F)$, and therefore $P\prec_M A\rtimes$Norm$(F)$. 
\hfill$\blacksquare$

We will also need the following lemma which is in the spirit of \cite[Lemma 4.2]{Va08}.

\begin{lemma}\label{commutant}\label{smaller}
Let $Q\subset qA_0^Fq$ be a von Neumann subalgebra for a subset $F\subset I$. Let $q'\in Q'\cap qAq$ be a non-zero projection such that $Qq'\nprec_{A}A_0^G$ for any proper subset $G\subset F$.

Then $q'(Q'\cap qMq)q'\subset q'(A\rtimes \text{Norm}(F))q'$.
 
In particular, if $Q\subset qA_0q$ is a diffuse subalgebra, then $Q'\cap qMq\subset q(A\rtimes\Sigma)q$.
\end{lemma}

{\it Proof.} Let $x\in q'(Q'\cap qMq)q'$ and consider $x=\sum_{g\in\Gamma}x_gu_g$ its Fourier decomposition. Take $g\in\Gamma\setminus \text{Norm}(F)$ and let us show that $x_g=0.$
Note that $Qq'\nprec_A A_0^{gF}$. Indeed, by assuming the contrary,
we can apply \cite[Lemma 2.4(3)]{DHI16} and obtain a non-zero projection $q_1'\in Q'\cap qMq$ with $q_1'\leq q'$ such that $Qq_1'\prec^s_A A_0^{gF}$. Since $Qq_1'\prec_A^s A_0^F $ and $A_0^F$ and $A_0^{gF}$ are in a commuting square position inside $A$, we can use
\cite[Lemma 2.8(2)]{DHI16} to derive that $Qq'\nprec_{A} A_0^{gF\cap F}.$ This is a contradiction since $gF\cap F$ is a proper subset of $F$. Since $Qq'\nprec_A A_0^{gF}$, there exists a sequence of unitaries $(u_n)_n\subset \mathcal U(Q)$ such that $\|E_{A_0^{gF}}(xq'u_n)\|_2\to 0$ for every $x\in A.$ Since $x\in Q'\cap qMq$, we deduce that $u_nx_g=x_g\sigma_{g}(u_n)$, for any $n\ge 1.$ Using that $\sigma_g(u_n)\in A_0^{gF}$ and $x_g=x_g\sigma_g(q)=x_gq'$, we derive that
$$
\|E_{A_0^{gF}}(x_g^*x_g)\|_2=\|E_{A_0^{gF}}(x_g^*x_g\sigma_g(u_n))\|_2=
\|E_{A_0^{gF}}(x_g^*u_nx_g)\|_2\to 0,
$$
as $n$ goes to infinity. This shows that $E_{A_0^{gF}}(x_g^*x_g)=0$, hence $x_g=0$. This proves the lemma.
\hfill$\blacksquare$

\subsection{Weak mixing techniques}
Following \cite[Lemma 2.4]{PV06} we say that an action $\Gamma\car I$ of a group $\Gamma$ on a set $I$ is {\it weakly mixing} if for every $A,B\subset I$ finite, there exists $g\in\Gamma$ such that $gA\cap B=\emptyset$. 

\begin{lemma}\label{wm1}
Let $\Sigma<\Gamma$ be countable groups satisfying $[\Sigma: g\Sigma g^{-1}\cap \Sigma]=\infty$, for any $g\in\Gamma\setminus\Sigma$. Denote $I=\Gamma/\Sigma.$

Then the translation action $\Sigma\car I\setminus\{\Sigma\}$ is weakly mixing.
\end{lemma}

{\it Proof.} By using \cite[Lemma 2.4]{PV06}, it is enough to show that every orbit of $\Sigma\car I\setminus\{\Sigma\}$ is infinite. Let $g\Sigma\in I$ with $g\notin\Sigma$ and note that the stabilizer of $g\Sigma$ equals $g\Sigma g^{-1}\cap\Sigma$. Since $[\Sigma: g\Sigma g^{-1}\cap \Sigma]=\infty$, if follows that the orbit of $g\Sigma$ is infinite, hence proving the lemma.
\hfill$\blacksquare$




A trace preserving action $\Gamma\overset\sigma\car (A,\tau)$ is called {\it weakly mixing} if for any $a_1,\dots,a_n\in A$ and $\epsilon>0$, there exists $g\in\Gamma$ such that $|\tau(a_i\sigma_g(a_j))-\tau(a_i)\tau(a_j) |<\epsilon$, for all $1\leq i,j\leq n.$

The proof of the following lemma is standard and we leave it to the reader.

\begin{lemma}\label{wm2}
Let $\Sigma<\Gamma$ be countable groups. Let $\Sigma\car A_0$ be a trace preserving action and denote by $\Gamma\overset{\sigma}{\car} A_0^{I}$ the associated coinduced action. Let $I_0=I\setminus\{\Sigma\}$ and note that $\sigma_{|\Sigma}$ is isomorphic to the diagonal product action $\Sigma\car A_0\bar\otimes A_0^{I_0}$. Assume the translation action $\Sigma\car I_0$ is weakly mixing. 

Then the restricted action $\Sigma \car A_0^{I_0}$ is weakly mixing. 
\end{lemma}

For a p.m.p. weakly mixing action $\Gamma\car (X,\mu)$, the quasi-normalizer of $L(\Gamma)$ inside $L^\infty(X)\rtimes\Gamma$ equals $L(\Gamma)$
\cite[Proposition 6.10]{Io08a}. The following lemma extends this result to diagonal product actions and its proof is contained in the proof of \cite[Theorem 3.1]{Po03}.

\begin{lemma}[\!\!\cite{Po03}]\label{L:qn}\label{wm3}
Let $\Gamma\car A$ and $\Gamma\car B$ be trace preserving actions and denote $M=(A\bar\otimes B)\rtimes\Gamma$ and $N=A\rtimes\Gamma$. Assume $\Gamma\car B$ is weakly mixing.

If there exist elements $x,x_1,\dots,x_n\in M$ such that $Nx\subset \sum_{i=1}^nx_iN$, then $x\in N$.
In particular, the quasi-normalizer of $N$ inside $M$ equals $N$.

\end{lemma}





\section{Rigidity results for coinduced actions}

We start this section by recalling the free product deformation for coinduced actions which was originally introduced for general Bernoulli actions by A. Ioana \cite{Io06a}.

\begin{definition}\label{def:free}
Let $\Gamma\car A_0^I$ be the coinduced action of a trace preserving action $\Sigma\overset{\sigma_0}\car A_0$ and let $M=A_0^I\rtimes\Gamma$. Consider the free product $\tilde A_0:=A_0*L(\mathbb Z)$ with respect to the natural traces and extend canonically $\Sigma\overset{\sigma_0}\car A_0$ to an action on $\tilde A_0$, still denoted by $\tilde\sigma_0$. More precisely, $(\tilde\sigma_0)_g(a_1b_1\dots a_nb_n)=(\sigma_0)_g(a_1)b_1\dots (\sigma_0)_g(a_n)b_n$, for all $g\in\Sigma$, $a_1,\dots, a_n\in A_0$ and $b_1,\dots,b_n\in L(\mathbb Z)$.
We denote by $\Gamma\car \tilde A_0^I$ the coinduced action associated to $\Sigma_0\car\tilde A_0$ and let $\tilde M=\tilde A_0^I\rtimes\Gamma.$ 

Define the self-adjoint $h\in L(\mathbb Z)$ with spectrum $[-\pi,\pi]$ such that exp$(ih)$ equals the canonical generating unitary $v\in L(\mathbb Z)$. Let $v_t=$ext$(ith)$ and denote by $\alpha_t^0\in $Aut$(\tilde A_0)$ the inner automorphism given by $\alpha_t^0(x)=v_txv_t^*$, for all $t\in\mathbb R$ and $x\in \tilde A_0$.
Denote by $\alpha_t\in$Aut$(\tilde A_0^{I})$ the automorphism given by $\alpha_t=\otimes_{i\in I}\alpha_t^0$. Since $\alpha_t$ is $\Gamma$-equivariant, we can extend it in a natural way to an automorphism of $\tilde M$ by letting $\alpha_t(u_g)=u_g$, for any $g\in\Gamma$. 
\end{definition}

For the next theorem, we assume the setting of Definition \ref{def:free}. The result goes back to Popa's spectral gap argument \cite[Lemma 5.1]{Po06} and it is similar to \cite[Corollary 4.3]{IPV10} and \cite[Theorem 3.1]{BV12}. We omit its proof since it is identical to the proof of \cite[Theorem 3.1]{Dr17}.

\begin{theorem}\label{spectral gap}
Assume $A_0$ is amenable and $g\Sigma g^{-1}\cap\Sigma$ is amenable for any $g\in\Gamma\setminus\Sigma$. Let $N$ be a tracial von Neumann algebra and let $Q\subset q(M\bar\otimes N)q$ be a von Neumann subalgebra such that
 $Q'\cap qMq$ is strongly non-amenable relative to $1\otimes N$.
 
 Then $\alpha_t\otimes id\to id$ uniformly as $t\to 0$ on $(Q)_1$.

\end{theorem}

We also need the following variant of \cite[Theorem 4.2]{IPV10} adapted to restricted coinduced actions. The result goes back to \cite[Theorem 4.1]{Po03} and \cite[Theorem 3.6]{Io06a}, and  since our proof is similar, we are rather brief.

\begin{theorem}\label{rigid2}
Let $\Sigma<\Gamma$ be countable groups. Let $\Sigma\car A_0$ be a trace preserving action and let $\Gamma\car A:=A_0^I$ be the associated coinduced action. Let $F_0\subset I$ be a finite, possibly empty, subset and denote $M=A\rtimes$Norm$(F_0)$.

Let $Q\subset qMq$ be a von Neumann subalgebra such that $\alpha_t\to id$ uniformly as $t\to 0$ on $(Q)_1$ and denote by $P$ the quasi-normalizer of $Q$ inside $qMq$. Then one of the following holds:

\begin{itemize}

\item $P\prec_M A\rtimes Stab(F_0,i)$ for some $i\in I\setminus F_0$.

\item $Q\prec_M A_0^{F}\rtimes Stab(F)$ for some finite subset $F_0\subset F\subset I.$ 

\end{itemize}

Moreover, if there exists $k\in\mathbb N$ such that $Stab(F)$ is finite for any subset $F_0\subset F\subset I$ with $|F|\ge k$, then of the following holds: 

\begin{itemize}

\item $Q\prec_M A_0^{F_0}$.

\item $P\prec_M A\rtimes$Stab$(F_0,i)$, for some $i\in I\setminus F_0$.

\item There exists a non-zero partial isometry $v\in Mq$ with $v^*v\in P$ and $vPv^*\subset A_0^{F_0}\rtimes$Norm$(F_0)$. In this case, we can assume that $v^*v=q$, provided $A_0^{F_0}\rtimes$Norm$(F_0)$ is a II$_1$ factor.

\end{itemize}

\end{theorem}

Before proceeding to the proof, we clarify the notation we used and then make a remark. We denoted by $\alpha_t$ the free product deformation associated to the coinduced action $\Gamma\car A$ as in Definition \ref{def:free}. For $i\in I$, we denoted by Stab$(F_0,i)$ the stabilizer of $F_0\cup \{i\}$ with respect to the action $\Gamma\car I$. We also used the convention that Norm$(\emptyset):=\Gamma$.

\begin{remark} We highlight here how Theorem \ref{rigid2} will be used for proving our main results.

\begin{itemize}

\item Theorem \ref{rigid2} will be applied in the cases $F_0=\{\emptyset\}$ and $F_0=\{\Sigma\}$. Note that these cases correspond to the classification of all ``rigid'' subalgebras of $A\rtimes\Gamma$ and $A\rtimes\Sigma$, respectively.

\item In both of these two cases, the assumption of the moreover part is satisfied by taking $k=2$ if we assume that
$g\Sigma g^{-1}\cap \Sigma$ is finite for any $g\in\Gamma\setminus\Sigma$. Indeed, note that if $F=\{g\Sigma,h\Sigma \}$ for some $g,h \in\Gamma$ such that $h^{-1}g \notin\Sigma$, then Stab$(F)=g\Sigma g^{-1}\cap h\Sigma h^{-1}$ is finite.

\item If $F_0=\{\emptyset\}$, the moreover part of Theorem \ref{rigid2} is showing that if $Q$ is diffuse, then $P\prec_{A\rtimes\Gamma} A\rtimes\Sigma$ or $P\prec_{A\rtimes\Gamma} L(\Gamma)$. On the other hand, if $F_0=\{\Sigma\}$, then the moreover part is showing that $Q\prec_{A\rtimes\Sigma} A_0$ or $P\prec_{A\rtimes\Sigma} A$ or $P\prec_{A\rtimes\Sigma} A_0\rtimes\Sigma$. 
\end{itemize}

\end{remark}


{\it Proof.} We define a free product deformation on $M$ as follows. Denote $I_0=I\setminus F_0$. Consider the free product $B_0=A_0*L(\mathbb Z)$ with respect to the natural traces. Note that the action $\Sigma\car A_0$ extends naturally (as in Definition \ref{def:free}) to a trace preserving action $\Sigma\car B_0$ and denote by $\Gamma\car B_0^I$ the associated coinduced action. 
It is clear that $A_0^{I_0}\subset B_0^{I_0}$ and that $\Sigma\car B_0^{I_0}$ extends $\Sigma\car A_0^{I_0}$, hence we have the natural inclusions $M\subset \tilde M:=(A_0^{F_0}\bar\otimes B_0^{I_0})\rtimes$Norm$(F_0)  \subset B_0^I\rtimes$Norm$(F_0)$.

Define the self-adjoint $h\in L(\mathbb Z)$ with spectrum $[-\pi,\pi] $ such that exp$(ih)$ equals the canonical generating unitary $v\in L(\mathbb Z)$. Let $v_t=$ext$(ith)$ and denote by $\alpha_t^0\in $Aut$(B_0)$ the inner automorphism given by $\alpha_t^0(x)=v_txv_t^*$, for all $t\in\mathbb R$ and $x\in B_0$.
Denote by $\alpha^1_t\in$Aut$(B_0^{I_0})$ the automorphism given by $\alpha^1_t=\otimes_{i\in I_0}\alpha_t^0$. Note that we can extend in a natural way $\alpha^1_t$ to an automorphism of $\tilde M$ by letting $\alpha^1_t(x)=x$, for any $x\in A_0^{F_0}$ and $\alpha^1_t(u_g)=u_g$, for any $g\in $ Norm$(F_0)$.

Assume now that $P\nprec_M A\rtimes \text{Stab}{(F_0,i)}$, for any $i\in I_0$. Since the deformation $\alpha_t$ equals to $(\otimes_{F_0}\alpha_t^0)\otimes\alpha_t^1$, by using the triangle inequality we derive that $\alpha^1_t\to id$ uniformly as $t\to 0$ on $(Q)_1$. Indeed, note first that for any $x\in Q$ we have
\begin{equation}\label{ti}
\|\alpha_t^1(x)-x\|_2\leq \|((\otimes_{F_0}\alpha_t^0)\otimes id) (\alpha_t^1(x))-\alpha_t^1(x)\|_2+\|\alpha_t(x)-x\|_2
\end{equation}
We have that $\alpha_t^0\to id$ uniformly as $t\to 0$ on $(A_0*L(\mathbb Z))_1$ by using, for instance, \cite[Lemma 3.6]{Io06b}. Since $F_0$ is a finite set, we derive that $(\otimes_{F_0}\alpha_t^0)\otimes id\to id$ uniformly as $t\to 0$ on $(\tilde M)_1$. Therefore, the inequality \eqref{ti} implies that $\alpha^1_t\to id$ uniformly as $t\to 0$ on $(Q)_1$.
Hence, we can find $t=2^{-n}$, for some positive integer $n$ and a non-zero partial isometry $v_t\in \tilde M$ such that $xv_t=v_t\alpha^1_t(x)$, for all $x\in Q.$ This proves that $v_t$ is a $Q$-$\alpha^1_t(Q)$-finite element. Here, we say that an element $x\in \mathcal M$ is $Q_1$-$Q_2$-finite for some given von Neumann algebras $Q_1,Q_2\subset \mathcal M$ if there exist $x_1,\dots,x_n,y_1,\dots,y_n\in\mathcal M$ such that 
$
xQ_2\subset \sum_{i=1}^n Q_1x_i \text{  and  }
Q_1x\subset \sum_{i=1}^n y_jQ_2.
$ 

We now show that there exists a non-zero $Q$-$\alpha^1_{2t}(Q)$-finite element $v_{2t}\in \tilde M$. Denote by $\mathcal D$ the set of all  $Q$-$Q$-finite elements in $q M q$ and note that  $\alpha^1_t(\beta(v_{t}^*)dv_{t})$ is a $Q$-$\alpha^1_{2t}(Q)$-finite element in $q\tilde M q$, for any $d\in\mathcal D$. Therefore, we just need to find $d\in\mathcal D$ satisfying $\beta(v_{t}^*)dv_{t}\neq 0.$ Suppose that $\beta(v_{t}^*)dv_{t}=0$, for any $d\in\mathcal D$. We denote by $q_1\in q\tilde Mq$ the unique projection onto the $\|\cdot\|_2$-closed linear span of $\mathcal Dv_t \tilde M$. Hence, $q_1\in P'\cap q\tilde Mq$ and $\beta(q_1)q_1=0$. Lemma \ref{quasi}(1) below shows that $q_1\in M$ and therefore $\beta(q_1)=q_1$. This implies that $q_1=0$, contradiction. Thus, there exists a non-zero $Q$-$\alpha^1_{2t}(Q)$-finite element $v_{2t}\in \tilde M$. Continuing inductively, we find a non-zero $Q$-$\alpha^1_{1}(Q)$-finite element $v_{1}\in \tilde M$, which shows that $\alpha_1(Q)\prec_{\tilde M}M.$

We continue by showing that there exists a finite subset $F_0\subset F\subset I$ such that $Q\prec_M A_0^{F}\rtimes \text{Stab}(F)$. We assume the contrary. Then, by applying an idea that goes back to \cite[proof of Theorem 4.3]{IPP05}, we obtain that there exists a sequence of unitaries $(u_n)_n\subset\mathcal U(Q)$ such that $\|E_{A_0^{F}\rtimes \text{Stab}(F)}(xu_ny)\|_2\to 0$, for all $x,y\in M$ and finite subsets $F_0\subset F\subset I$. 
By proceeding as in the proof of \cite[Theorem 4.2]{IPV10}, it follows that 
$\|E_M(x\alpha_1(v_n)y)\|_2\to 0$, for all $x,y\in M$, proving that $\alpha_1(Q)\nprec_{\tilde M}M,$ contradiction.

For proving the moreover part, we assume that $Q\nprec_M A_0^{F_0}$. Since $P\nprec_M A\rtimes \text{Stab}{(F_0,i)}$, for any $i\in I_0$, we need to show that $P\prec_M A_0^{F_0}\rtimes$Norm$(F_0)$. Note that the previous paragraph shows that there exists a finite subset $F_0\subset F\subset I$ such that $Q\prec_M A_0^{F}\rtimes \text{Stab}(F)$. We can assume that $F$ is minimal with this property. Therefore, we notice that $Q\nprec_M A_0^{F}$. Indeed, if $Q\prec_M A_0^{F}$, note that we can apply Lemma \ref{normalizer} and deduce that $P\prec A\rtimes$Norm$(F).$ Since $F$ is finite, we have that Stab$(F)$ is a finite index subgroup of Norm$(F)$, which implies that $P\prec A\rtimes$Stab$(F).$ Since $Q\nprec_M A_0^{F_0}$, we have that $F_0$ is strictly contained in $F$. Hence, there exists $i\in I_0$ such that $P\prec A\rtimes$Stab$(F_0,i),$ contradiction. Hence, $Q\nprec_M A_0^{F}$.

Note that there exists a finite subset $G\supset F$ such that $Q\prec_M A_0^F\rtimes$Stab$(G)$ and $Q\nprec_M A_0^F\rtimes$Stab$(G')$, for any $G'$ strictly larger than $G$.
Therefore, there exist some projections $q_0\in Q, r\in A_0^{F}\rtimes$Stab$(G)$, a non-zero partial isometry $v\in rMq_0$ and a $*$-homomorphism $\phi: q_0Qq_0\to r(A_0^{F}\rtimes$Stab$(G))r$ such that $\phi(x)v=vx$, for any $x\in q_0Qq_0$ and such that 
\begin{equation}\label{aa11}
\phi(q_0Qq_0)\nprec_{A_0^{F}\rtimes \text{Stab}(G)} A_0^{F}\rtimes\text{Stab}(G'), \text{  for any } G' \text{  strictly larger than }G.
\end{equation}

Next, we show that $G=F_0$. Assume this is not the case. Relation \eqref{aa11} implies that
$$
\phi(q_0Qq_0)\nprec_{A\rtimes \text{Stab}(G)} A\rtimes\text{Stab}(G'), \text{  for any } G' \text{  strictly larger than }G.
$$
We can use Lemma \ref{quasi}(2) bellow and deduce that $P\prec_M A\rtimes$Stab$(G)$. Since $F_0\subsetneq G$, we derive that there exists $i\in I_0$ such that $P\prec_M A\rtimes$Stab$(F_0,i)$, contradiction. Therefore, $G=F_0$, and hence, $F=F_0$. 

Finally, note that $\phi(q_0Qq_0)\subset r(A_0^{F_0}\rtimes$Stab$(F_0))r$. Since $v^*v$ commutes with $q_0Qq_0$, it belongs to $P$.
By applying Lemma \ref{quasi}(1) and using \eqref{aa11} we obtain that $vPv^*\subset A_0^{F_0}\rtimes$Norm$(F_0)$. Assume now that $A_0^{F_0}\rtimes$Norm$(F_0)$ is a II$_1$ factor. Take partial isometries $v_1,\dots,v_k\in P$ such that $v_i^*v_i\leq v^*v$ and $\sum_{i=1}^k v_iv_i^*$ is a central projection in $P$. Since $A_0^{F_0}\rtimes$Norm$(F_0)$ is a II$_1$ factor, there exist partial isometries $w_1,\dots,w_k\in A_0^{F_0}\rtimes$Norm$(F_0)$ such that $w_iw_i^*=vv_i^*v_iv^*$ and $w_iw_j^*=0$, for all $i\neq j$. Define $w=\sum_{i=1}^k v_iv^*w_i\in qM$ and note that $w$ is a partial isometry satisfying $ww^*\in \mathcal Z(P)$ and $w^*Pw\subset A_0^{F_0}\rtimes$Norm$(F_0)$. 
Since $A_0^{F_0}\rtimes$Norm$(F_0)$ is II$_1$ factor, a standard maximality argument (see, for instance, the second paragraph of the proof of \cite[Corollary 4.3]{IPV10}) allows us to construct a partial isometry 
$v_0\in Mq$ with $v_0^*v_0=q$ that satisfies $v_0Pv_0^*\subset A_0^{F_0}\rtimes$Norm$(F_0)$
\hfill$\blacksquare$





We continue with the following lemma which was needed in the proof of Theorem \ref{rigid2}. We omit the proof since it is proved exactly in the same way as \cite[Lemma 4.1]{IPV10}.

\begin{lemma}\label{quasi}
Let $\Sigma<\Gamma$ be countable groups. Let $A_0$, $B_0\subset C_0$ be tracial von Neumann algebras and let $\Sigma\car A_0$ and $\Sigma\car C_0$ be trace preserving actions such that $B_0$ is $\Sigma$-invariant. Denote by $\Gamma\car A_0^I$, $\Gamma\car B_0^I$ and $\Gamma\car C_0^I$ the coinduced actions associated to $\Sigma\car A_0$, $\Sigma\car B_0$ and $\Sigma\car C_0$, respectively.

Let $F_0\subset I$ be a finite, possibly empty, subset and let $I_0=I\setminus F_0$. We consider the diagonal product actions Norm$(F_0)\car A_0^{F_0}\bar\otimes B_0^{I_0}$ and Norm$(F_0)\car A_0^{F_0}\bar\otimes C_0^{I_0}$ and denote by $M=(A_0^{F_0}\bar\otimes B_0^{I_0})\rtimes$Norm$(F_0)$ and $\tilde M=(A_0^{F_0}\bar\otimes C_0^{I_0})\rtimes$Norm$(F_0)$. Note that $M\subset\tilde M$ and let $A=A_0^{F_0}\bar\otimes B_0^{I_0}.$

\begin{enumerate}

\item If $P\subset pMp$ is a von Neumann subalgebra such that $P\nprec_M A\rtimes Stab(F_0,i)$ for any $i\in I_0$, then the quasi-normalizer of $P$ inside $p\tilde Mp$ is contained in $pMp$.

\item If $G\supset F_0$ is a finite set and $Q\subset q(A\rtimes\text{Stab}( G))q$ is a von Neumann subalgebras such that $Q\nprec_{A\rtimes\text{Stab}(G)}A\rtimes\text{Stab}(G')$, for any strictly larger subset $G'\supset G$, then the quasi-normalizer of $Q$ inside $qMq$ is contained in 
$q(A\rtimes\text{Norm}( G))q$.

\end{enumerate}

\end{lemma}

Finally, we prove the following lemma that explores the rigidity assumption of non-amenable groups which are measure equivalent in the sense of Gromov \cite{Gr91} to products of infinite groups. Since the proof is standard, we will only provide a sketch. We will use the following characterization of measure equivalency due to Furman \cite{Fu99}.
\begin{definition}\label{def:ME}
Two countable groups are {\it measure equivalent} if they admit free ergodic stable orbit equivalent p.m.p. actions. Natural examples of measure equivalent groups are provided by pairs of lattices in an unimodular locally compact second countable group.
\end{definition}

\begin{lemma}\label{L:ME}
Let $\Sigma<\Gamma$ be countable groups such that $g\Sigma g^{-1}\cap \Sigma$ is amenable for any $g\in\Gamma\setminus\Sigma$. Let $\Sigma\car (X_0,\mu_0)$ be a p.m.p. action and let $\Gamma\car (X,\mu)$ be the corresponding coinduced action. Denote $M=L^\infty(X)\rtimes\Gamma$ and let $\alpha_t$ be the associated free product deformation as in Definition \ref{def:free}.

Let $Q\subset qMq$ be a von Neumann subalgebra such that $Q=L(\Lambda)$, where $\Lambda$ is a group that is measure equivalent to a product of two non-amenable groups.

Then $\alpha_t\to id$ uniformly on $(Q)_1$ 
\end{lemma}

{\it Proof.} Assume for simplicity that there exist free ergodic p.m.p. actions of $\Lambda$ and $\Lambda_1\times\Lambda_2$, with $\Lambda_1$ and $\Lambda_2$ non-amenable groups, on a measurable space $(Y,\nu)$ whose orbits are equal, almost everywhere.
Denote $N=L^\infty(Y)\rtimes\Lambda=L^\infty(Y)\rtimes(\Lambda_1\times
\Lambda_2)$ and let $\Delta:N\to N\bar\otimes L(\Lambda)$ be the $*$-homomorphism defined by $\Delta(bv_\lambda)=bv_\lambda\otimes v_\lambda$, for all $b\in L^\infty(Y)$ and $\lambda\in\Lambda$. 

Note that $\Delta(L(\Lambda_i))$ is strongly non-amenable relative to $N\otimes 1$ by \cite[Lemma 10.2(5)]{IPV10} for any $1\leq i\leq 2$. Using Theorem \ref{spectral gap}, we derive that $id\otimes \alpha_t\to id$ uniformly on $(\Delta(L(\Lambda_i)))_1$, for any $1\leq i\leq 2$, and hence, on $(\Delta(N))_1$. This clearly implies that $\alpha_t\to id$ uniformly on $(Q)_1$.
\hfill$\blacksquare$

\section{Proofs of Theorem \ref{A:solid}, Proposition \ref{strong ergodicity} and Theorem \ref{B:solid}}

The proofs of Theorem \ref{A:solid} and Theorem \ref{B:solid} need the following useful characterizations for the two notions of solidity \cite[Proposition 6]{CI08}.

\begin{proposition}[\!\!\cite{CI08}]\label{equivalence}
(i) A measure preserving equivalence relation $\mathcal R$ on a probability space $(X,\mu)$ is solidly ergodic if and only if $Q'\cap L(\mathcal R)$ is amenable for any diffuse subalgebra $Q\subset L^\infty(X)$.

(ii) A tracial von Neumann algebra $M$ is solid if and only if $Q'\cap M$ is amenable for any diffuse subalgebra $Q\subset M$.

\end{proposition}

\subsection{\bf Proof of Theorem \ref{A:solid}}

Let $I=\Gamma/\Sigma$ and define the von Neumann algebras $A_0=L^\infty(X_0), A=L^\infty(X)$ and $ M=A\rtimes\Gamma$. 
Let $Q\subset A$ be a diffuse von Neumann subalgebra. By Proposition \ref{equivalence}, we only have to show that $Q'\cap M$ is amenable.
Assume by contradiction that $Q'\cap M$ is non-amenable.
A maximality argument combined with \cite[Lemma 2.6]{DHI16} implies that
 there exists a non-zero projection $z\in \mathcal Z(Q'\cap M)$ such that $(Q'\cap M)z$ is strongly non-amenable. By applying Theorem \ref{spectral gap}, we obtain that $\alpha_t\to id$ uniformly as $t\to 0$ on $(Qz)_1$. Let $a\in A$ be the support projection of $E_A(z)$ and remark that $a_1z\neq 0$, for any non-zero projection $a_1\in Aa.$ Therefore, since $(Q'\cap M)z$ is properly non-amenable and $z\in A'\cap M$, it follows that 
\begin{equation}\label{a1}
a_1(Q'\cap M)a_1 \text{ is non-amenable for any non-zero projection } a_1\in Aa.
\end{equation}

By applying Theorem \ref{rigid2}, it follows that there exists a finite set $F\subset I$ such that $Qz\prec_M A_0^F\rtimes$Norm$(F)$. Hence, by \cite[Lemma 2.4(3)]{DHI16} there exists a non-zero projection $z_1\in (Q'\cap M)z$ such that $Qz_1\prec_M^s A_0^F\rtimes$Norm$(F)$. Since $Q\subset A$, \cite[Remark 2.2]{DHI16} implies that $Qz\prec_M^s A$. Since $A$ is regular in $M$ and $A$ and $A_0^F\rtimes$Norm$(F)$ are in a commuting square position, we can apply \cite[Lemma 2.8(2)]{DHI16} and derive that $Qz_1\prec_M^s A_0^F.$ Hence, $Qa\prec_M A_0^F$ since $z_1\leq a$. We can apply Lemma \ref{inside} and deduce that $Qa\prec_A A_0^{F}$ by replacing, eventually, $F$ by another finite set.
Moreover, we can assume that $F\subset I$ is a minimal subset with the property that there exists a non-zero projection $a_0\in A$ such that $Qa_0\prec_A A_0^{F}.$

Since $A$ is abelian and $Qa\prec_A A_0^F$, there exist non-zero projections $a_1\in Qa, r\in A_0^F, a_2\in Aa_1r$ and a one-to-one $*$-homomorphism $\varphi: Qa_1\to A_0^Fr$ such that  $\varphi(x)a_2=xa_2$, for every $x\in Qa_1.$
Note that $R:=\varphi(Qa_1)\subset A_0^Fr$ satisfies
\begin{equation}\label{a2} 
Ra_2=Qa_2 \text{ and } a_2(R'\cap rMr)a_2=a_2(Q'\cap M)a_2.
\end{equation}

We may assume that 
$ Ra_2\nprec_{A} A_0^G$ for any proper subset $G\subset F$.
Indeed, if there exists such a finite set $G$, then we would have that $Qa_2\prec_A A_0^G$, which contradicts the minimality of $F$. Hence, by applying Lemma \ref{commutant}, we deduce that $a_2(R'\cap rMr)a_2=a_2(R'\cap r(A\rtimes\text{Norm}(F))r)a_2.$

Next, we use relations \eqref{a1} and \eqref{a2} to derive that $a_2(R'\cap rMr)a_2$ is non-amenable, and hence, $R'\cap r(A\rtimes\text{Norm}(F))r$ is non-amenable. Since $g\Sigma g^{-1}\cap \Sigma$ is amenable for any $g\in\Gamma\setminus\Sigma$, we obtain that $F$ has only one element. Let $F=\{g\Sigma\}$ for some $g\in\Gamma$.
Note that Norm$(F)=g\Sigma g^{-1}$ and $R\subset A_0^{g\Sigma}r.$ It follows that $R'\cap r(A\rtimes g\Sigma g^{-1})r$ is non-amenable.
Since the action $g\Sigma g^{-1}\car I$ is fixing the element $g\Sigma$, we have that $g\Sigma g^{-1}\car A$ is isomorphic to the diagonal action $g\Sigma g^{-1}\car A_0^{g\Sigma}\bar\otimes A_0^{I\setminus\{g\Sigma\}}$. Therefore, we can apply Proposition \ref{splits} and conclude that $R'\cap r(A_0^{g\Sigma}\rtimes g\Sigma g^{-1})r$ is non-amenable. Finally, note that $\tilde R:=R\oplus A_0^{g\Sigma}(1-r)$ is a diffuse von Neumann subalgebra of $A_0^{g\Sigma}$ for which its relative commutant $\tilde R'\cap ( A_0^{g\Sigma}\rtimes g\Sigma g^{-1})$ is non-amenable. This contradicts the solid ergodicity of $\mathcal R_{\Sigma\car X_0}$, which proves the theorem.
\hfill$\blacksquare$


\begin{remark}
In the proof of Theorem \ref{A:solid}, there is no need to consider the support projection of $E_A(z)$ if the coinduced action $\Gamma\car (X,\mu)$ is free since $z\in\mathcal Z(Q'\cap M)\subset A'\cap M=A$ in this case. Note that the freeness condition is guaranteed if $\cap_{g\in\Gamma}g\Sigma g^{-1}$ is trivial  or if $\Sigma\car (X_0,\mu_0)$ is free, see \cite[Lemma 2.1]{Io06c} and \cite[Lemma 5.3]{Dr15}.
\end{remark}



\subsection{\bf Proof of Proposition \ref{strong ergodicity}}

Let $\{e_i\}_{i\ge 0}$ be an orthonormal basis for $L^2(X_0)$, with $e_0=1$. Thus $\{\otimes_{i\in F} e_{i}|F\subset \Gamma/\Sigma \text{ finite subset} \}$ is an orthonormal basis for $L^2(X)$.
Let $S$ be the countable family of finite non-empty subsets of $\Gamma/\Sigma$. Therefore, any $\xi\in L^2(X)\ominus\mathbb C 1$ admits a unique canonical decomposition $\xi=\sum_{F\in S} \xi_F$, where $\xi_F\in \otimes_F (L^2(X_0)\ominus\mathbb C1)$.
We define the map $\theta:L^2(X)\ominus\mathbb C1\to \ell^2(S)_{+}$ by letting $\theta(\xi)=\sum_{F\in S}\|\xi_F\|_2\delta_F$, for any $\xi\in L^2(X)\ominus\mathbb C1$, and notice that $\|\theta(\xi)\|_2=\|\xi\|_2$.

We assume by contradiction that $\Gamma\overset\sigma\car X$ does not have spectral gap. Thus, there exists a sequence of unit vectors $(\xi_n)_n\subset L^2(X)\ominus\mathbb C1$ satisfying $\|\sigma_g(\xi_n)-\xi_n\|_2\to 0$, for any $g\in\Delta$.
 Consider the left translation action $\Gamma\overset\rho\car S$ and note that the sequence $\theta(\xi_n)$ is a sequence of unit vectors that satisfies $\|\rho_g(\theta(\xi_n))-\theta(\xi_n)\|_2\to 0$, for any $g\in\Delta$. Indeed, this follows from the fact that
$$
\|\rho_g(\theta(\xi_n))-\theta(\xi_n)\|_2^2=
\sum_{F\in S} |\|\xi_{n,F}\|_2-\|\xi_{n,gF}\|_2|^2\leq \sum_{F\in S}\|\sigma_g(\xi_{n,F})-\xi_{n,gF}\|_2^2=\|\sigma_g(\xi_n)-\xi_n\|_2^2\to 0,
$$
for any $g\in\Delta$. We remark that there is a natural $\Gamma$-equivariant identification of $S\cup\{\emptyset\}$ with the infinite direct sum $(\mathbb Z/2\mathbb Z)^{(\Gamma/\Sigma)}$. This identification maps $\{\emptyset\}$ to the neutral element in $(\mathbb Z/2\mathbb Z)^{(\Gamma/\Sigma)}$.
Therefore, we can conclude that the generalized Bernoulli action $\Delta\car (\mathbb Z/2\mathbb Z)^{(\Gamma/\Sigma)}$ does not have spectral gap. By applying \cite[Theorem 1.2]{KT08}, it follows that $\Delta$ is amenable relative to $\Sigma$ inside $\Gamma$, which is a contradiction.
\hfill$\blacksquare$

We continue by showing a von Neumann algebraic version of Proposition \ref{strong ergodicity}. 

\begin{proposition}\label{vNa version}

Let $\Sigma<\Gamma$ be countable groups satisfying $g\Sigma g^{-1}\cap \Sigma$ is finite for any $g\in\Gamma\setminus\Sigma$. Let $\Sigma{\car} (X_0,\mu_0)$ be any p.m.p. action and let $\Gamma\car (X,\mu)$ be the associated coinduced action. Denote by $M=L^\infty(X)\rtimes\Gamma$ and let $Q\subset qMq$ be a von Neumann subalgebra with property Gamma. 

If $Q$ is not amenable relative to $L^\infty(X)\rtimes\Sigma$ inside $M$, then $Q\prec_M L(\Gamma)$. 

In particular, if $\Delta<\Gamma$ is a subgroup that is not amenable relative to $\Sigma$, then $L^\infty(X)\rtimes\Delta$ does not have property Gamma.

\end{proposition}
 
{\it Proof.}
Let $I=\Gamma/\Sigma$ and define the von Neumann algebras  $A=L^\infty(X)$ and $ M=A\rtimes\Gamma$. Let $\alpha_t$ be the free product deformation associated to the coinduced action $\Gamma\car A$ as in Definition \ref{def:free}.
Since $Q$ has property Gamma, \cite[Theorem 3.1]{HU15} implies that there exists a decreasing sequence of diffuse abelian von Neumann subalgebras $A_n$ of $Q$ such that $Q=\vee_{n\ge 1} (A_n'\cap Q)$. We assume that $Q$ is not amenable relative to $A\rtimes\Sigma$ and prove that $Q\prec_M L(\Gamma)$.


Since relative amenability is closed under inductive limits \cite[Lemma 2.7]{DHI16}, it follows that there exists an $n$ such that $A_n'\cap Q$ is not amenable relative to $A\rtimes\Sigma$ inside $M$. By using \cite[Lemma 2.6]{DHI16} combined with a maximality argument, we can find a non-zero projection $z\in \mathcal N_{M}(A_n'\cap Q)'\cap M$ such that $(A_n'\cap Q)z$ is strongly non-amenable relative to $A\rtimes\Sigma$. In particular, we can apply Theorem \ref{spectral gap} and derive that $\alpha_t\to id$ uniformly on $(A_nz)_1$. By applying Theorem \ref{rigid2} in the case $F_0=\{\Sigma\},$ we have (1) $(A_n'\cap Q)z\prec_M A\rtimes\Sigma$ or (2) there exists a partial isometry $v\in Mz$ such that $v^*v\in (A_n'\cap Q)z$ and $v(A_n'\cap Q)zv^*\subset L(\Gamma)$.

Assume (1) holds. By \cite[Lemma 2.4(3)]{DHI16}, there exists a non-zero projection $z_1\in (A_n'\cap Q)'\cap M$ with $z_1\leq z$ such that $(A_n'\cap Q)z_1\prec_M^s A\rtimes\Sigma.$ By \cite[Lemma 2.6(3)]{DHI16} we derive that $(A_n'\cap Q)z_1$ is amenable relative to $A\rtimes\Sigma$, which is a contradiction.

Hence, (2) holds. This means that there exists a partial isometry $v\in Mz$ such that $v^*v\in (A_n'\cap Q)z$ and $v(A_n'\cap Q)zv^*\subset L(\Gamma)$. Let $m\ge n$ and note that $z\in A_n'\cap Q\subset A_m'\cap Q$. Assume by contradiction that
 $v(A_mz)v^*\prec_{M} A\rtimes\Sigma$. Note that $v(A_mz)v^*\nprec_{M} A$ since $v(A_mz)v^*\subset L(\Gamma)$. By using Lemma \ref{mixing}(2), we get that $vz(A_m'\cap Q)zv^*\prec_{M} A\rtimes\Sigma$. By applying \cite[Lemma 2.6(3)]{DHI16} as before, we derive again a contradiction. Hence, $v(A_mz)v^*\nprec_{M} A\rtimes\Sigma$, which implies that $v(A_mz)v^*\nprec_{M} L(\Sigma)$. Now we apply \cite[Corollary 3.6(1)]{Dr17} and derive that $vz(A_m'\cap Q)zv^*\subset L(\Gamma)$, for any $m\ge n$. This shows that $zQz\prec_M L(\Gamma)$, which proves the main result.
 
The last part of the proposition follows by applying \cite[Lemma 2.5]{BV12}.
\hfill$\blacksquare$

\subsection{\bf Proof of Theorem \ref{B:solid}}
Let $I=\Gamma/\Sigma$ and define $A_0=L^\infty(X_0), A=L^\infty(X)$ and $ M=A\rtimes\Gamma$. Define the free product deformation associated to the coinduced action $\Gamma \car A$ as in Definition \ref{def:free}.
In light of Proposition \ref{equivalence}, note that $M$ is solid if and only if $Q'\cap  M$ is amenable for any diffuse abelian subalgebra $Q\subset  M$. This follows from the fact that any diffuse von Neumann algebra contains a diffuse abelian subalgebra.


Therefore, take a diffuse abelian subalgebra $Q\subset M$ and assume by contradiction that $Q'\cap M$ is non-amenable. Let $z\in\mathcal Z(Q'\cap M)$ be a non-zero projection such that $(Q'\cap M)z$ is strongly non-amenable. By applying Theorem \ref{spectral gap}, it follows that $\alpha_t\to id$ uniformly on $(Qz)_1$.  By using the moreover part of Theorem \ref{rigid2}, we obtain that $(Q'\cap M)z\prec_M A\rtimes\Sigma$ or $(Q'\cap M)z\prec_M L(\Gamma).$ By applying Lemma \ref{solid}, the second possibility implies that $L(\Gamma)$ is not solid, contradiction.  

Hence, $(Q'\cap M)z\prec_M A\rtimes\Sigma$. Thus, there exist projections $z_1\in Q'\cap M, p\in A\rtimes\Sigma$, with $z_1\leq z$, a non-zero partial isometry $v\in pMz_1$ and a one-to-one $*$-homomorphism $\phi: z_1(Q'\cap M)z_1\to p(A\rtimes\Sigma)p$ satisfying $\phi(x)v=vx$ for any $x\in z_1(Q'\cap M)z_1$. Moreover, we may assume that there exists a positive $c$ such that $E_{A\rtimes\Sigma}(vv^*)\ge cp$. 

We continue by showing that $\alpha_t\to id$ uniformly on $(\phi(Qz_1))_1$. Since $E_{\tilde A\rtimes\Sigma}(vv^*)=E_{ A\rtimes\Sigma}(vv^*)$, we deduce that
$$
\begin{array}{rcl}
\| \alpha_t(\phi(x))E_{A\rtimes\Sigma}(vv^*)-\phi(x)E_{A\rtimes\Sigma}(vv^*) \|_2 &\leq& \|\alpha_t(\phi(x))v-\phi(x)v\|_2\\
&\leq& 2\|\alpha_t(v)-v\|_2+\|\alpha_t(x)-x\|_2,
\end{array}
$$
for any $x\in (Qz_1)_1$ and $t\in\mathbb R$. Since $E_{A\rtimes\Sigma}(vv^*)\ge cp$, there exists $x_0\in A\rtimes\Sigma$ such that $E_{A\rtimes\Sigma}(vv^*)x_0=p$. Thus, we obtain that $\alpha_t\to id$ uniformly on $(\phi(Qz_1))_1$.

Denote $M_0=A\rtimes\Sigma$. Since $\phi(z_1(Q'\cap M)z_1)\subset pM_0p$ is non-amenable, we have that $P:= \phi (Qz_1)'\cap pM_0p$ is non-amenable.
By applying \cite[Lemma 2.6]{DHI16}, there exists a non-zero projection $p'\in \mathcal N_{pM_0p}(P)'\cap pM_0p$  such that $Pp'$ is strongly non-amenable in $M_0$. Since $\alpha_t\to id$ uniformly on $(\phi(Qz_1))_1$, we
apply Theorem \ref{rigid2} in the case $F_0=\{\Sigma\}$ and obtain that $\phi(Qz_1)p'\prec_{M_0} A_0$ or $Pp'\prec_{M_0} A$ or $Pp'\prec_{M_0} A_0\rtimes\Sigma$. If any of the last two possibilities holds, we can apply Lemma \ref{solid} and obtain that $A$ or $A_0\rtimes\Sigma$ is not solid, since $Pp'$ is properly non-amenable. Therefore, $\phi(Qz_1)p'\prec_{M_0} A_0$. Finally, Lemma \ref{soliderg} implies that $\mathcal R_{\Sigma\car X_0}$ is not solidly ergodic since $Pp'$ is properly non-amenable. This is a contradiction.
\hfill$\blacksquare$

\section{Proof of Theorem \ref{A:converse}}

We will need the following result for the proof of Theorem \ref{A:converse}, which is, for instance, a corollary of \cite[Appendix A]{Bo14}.

\begin{lemma}\label{mixing}
Let $\Sigma<\Gamma$ be countable groups satisfying $g\Sigma g^{-1}\cap \Sigma$ is finite for any $g\in\Gamma\setminus\Sigma$. Let $\Gamma\car A$ be a trace preserving action and denote $M=A\rtimes\Gamma$. Let $Q\subset pMp$ be a von Neumann subalgebra such that $Q\nprec_M A$ and denote by $P$ the quasi-normalizer of $Q$ in $pMp$. Then the following hold:

\begin{enumerate}

\item If $Q\subset A\rtimes\Sigma$, then $P\subset A\rtimes\Sigma.$

\item If $Q\prec_M A\rtimes\Sigma$, then $P\prec_M A\rtimes\Sigma$. Moreover, if $A\rtimes\Sigma$ is a factor and $Q\prec^s_M A\rtimes\Sigma$, then there exists a unitary $u\in M$ such that $uPu^*\subset A\rtimes\Sigma.$

\end{enumerate}
\end{lemma}

\subsection {Proof of Theorem \ref{A:converse}}

The following result implies in particular Theorem \ref{A:converse}.

\begin{theorem}\label{general}
For any $1\leq i\leq 2$, let $\Sigma_i<\Gamma_i$ be countable non-amenable groups such that
 $g\Sigma_i g^{-1}\cap\Sigma_i$ is finite for any $g\in\Gamma_i\setminus \Sigma_i$.  For any $1\leq i\leq 2$, let $\Sigma_i\car (X_i,\mu_i)$ be a free ergodic p.m.p. action that belongs to one of the following classes:

\begin{itemize}

\item $\Sigma_i\car (X_i,\mu_i)$ is a rigid or compact action; $\Sigma_i$ has property (T) or it is measure equivalent to a product of two non-amenable groups.

\item $\Sigma_i\car (X_i,\mu_i)$ is conjugate to a product action of two infinite groups.
\end{itemize}

If the coinduced actions $\Gamma_1\car (X_1,\mu_1)^{\Gamma_1/\Sigma_1}$ and $\Gamma_2 \car (X_2,\mu_2)^{\Gamma_2/\Sigma_2}$ are orbit equivalent, then $\Sigma_1\car (X_1,\mu_1)$ and $\Sigma_2\car (X_2,\mu_2)$ are orbit equivalent.
\end{theorem}

{\it Proof of Theorem \ref{general}.}
Denote $Y_1=X_1^{\Gamma_1/\Sigma_1}$ and $Y_2=X_2^{\Gamma_2/\Sigma_2}.$ The assumption implies by \cite{Si55,FM75} that there exists an isomorphism of von Neumann algebras between $L^\infty(Y_1)\rtimes\Gamma_1$ and $L^\infty(Y_2)\rtimes\Gamma_2$ which sends $L^\infty(Y_1)$ onto $L^\infty(Y_2)$. For ease of notation, assume that $M:=L^\infty(Y_1)\rtimes\Gamma_1=L^\infty(Y_2)\rtimes\Gamma_2$ and $L^\infty(Y_1)=L^\infty(Y_2)$ and let us show the following claim. 

{\bf Claim 1.} We have $L^\infty(Y_1)\rtimes\Sigma_1\prec_M L^\infty(Y_2)\rtimes\Sigma_2$ and $L^\infty(Y_2)\rtimes\Sigma_2\prec_M L^\infty(Y_1)\rtimes\Sigma_1.$

{\it Proof of Claim 1.} We will prove only the first intertwining relation.
Let $\alpha_t$ be the free product deformation as in Definition \ref{def:free} associated to the coinduced action $\Gamma_2\car Y_2$. 

Assume first that $\Sigma_1$ has property (T) or it is measure equivalent to a product of two non-amenable groups. Then by using Lemma \ref{L:ME}, we derive that $\alpha_t\to id$ uniformly on $(L(\Sigma_1))_1$. On one hand, if $\Sigma_1\car (X_1,\mu_1)$ is rigid, then 
$\alpha_t\to id$ uniformly on $(L^\infty(X_1)\rtimes\Sigma_1)_1$. By applying Theorem \ref{rigid2} in the case $F_0=\{\emptyset\}$, we obtain that $L^\infty(X_1)\rtimes \Sigma_1\prec_M L^\infty(Y_2)\rtimes\Sigma_2$ or $L^\infty(X_1)\rtimes \Sigma_1\prec_M L(\Gamma_2)$. The second option does not hold since $L^\infty(X_1)$ is diffuse. Hence, we obtain in particular that $L(\Sigma_1)\prec_M L^\infty(Y_2)\rtimes\Sigma_2$. By using \cite[Lemma 2.3]{BV12}, we obtain that $L^\infty(Y_1)\rtimes\Sigma_1\prec_M L^\infty(Y_2)\rtimes\Sigma_2$. On the other hand, if $\Sigma_1\car (X_1,\mu_1)$ is compact, we can apply \cite[Theorem 6.9]{Io08a} and derive that the quasi-normalizer of $L(\Sigma_1)$ contains $L^\infty(X_1)\rtimes\Sigma_1$. Once again, we use Theorem \ref{rigid2} and obtain that $L^\infty(X_1)\rtimes\Sigma_1\prec_M L^\infty(Y_2)\rtimes\Sigma_2$.


Assume now that $\Sigma_1\car X_1$ is isomorphic to $\Sigma_1^1\times\Sigma_1^2\car X_1^1\times X_1^2$, with $\Sigma^1_1$ non-amenable. \cite[Lemma 2.6(2)]{DHI16} combined with a maximality argument imply that we can find a non-zero projection $z\in (L^\infty(X_1)\rtimes\Sigma_1)'\cap M$ such that 
$z((L^\infty (X_1^2)\rtimes\Sigma_1^2 )'\cap M)z$ is strongly non-amenable.
By applying Theorem \ref{spectral gap} we have that
$\alpha_t\to id$ uniformly on $((L^\infty (X_1^2)\rtimes\Sigma_1^2 )z)_1$. By applying the moreover part of Theorem \ref{rigid2}, we obtain that $L^\infty(X_1)\rtimes\Sigma_1\prec_M L^\infty(Y_2)\rtimes\Sigma_2$ or $L^\infty(X_1)\rtimes\Sigma_1\prec_M L(\Gamma)$. The second option leads to the contradiction that $L^\infty(X_1)$ is not diffuse. Hence, $L(\Sigma_1)\prec_M L^\infty(Y_2)\rtimes\Sigma_2$.
By applying, for instance, \cite[Lemma 2.3]{BV12} we derive that $L^\infty(Y_1)\rtimes\Sigma_1\prec_M L^\infty(Y_2)\rtimes\Sigma_2$ since $\Sigma_1$ is normalizing $L^\infty(Y_1)$. 
\hfill$\square$

Note that Lemma \ref{wm1} and Lemma \ref{wm2} imply that
$\Sigma_1\car Y_1$ and $\Sigma_2\car Y_2$ are ergodic actions. Note also that the action $\Gamma_1\car Y_1$ is free by \cite[Lemma 2.1]{Io06c}.
Hence, the relative commutant $(L^\infty(Y_1)\rtimes\Sigma_1)'\cap M=\mathbb C1.$ Therefore, Claim 1 implies that $L^\infty(Y_1)\rtimes\Sigma_1\prec^s_M L^\infty(Y_2)\rtimes\Sigma_2$. By applying  Lemma \ref{mixing}(2), we derive that there exists a unitary $w_1\in M$ such that $w_1(L^\infty(Y_1)\rtimes\Sigma_1)w_1^*\subset L^\infty(Y_2)\rtimes\Sigma_2$. In a similar way, we can show that there exists a unitary $w_2\in M$ such that $w_2(L^\infty(Y_2)\rtimes\Sigma_2)w_2^*\subset L^\infty(Y_1)\rtimes\Sigma_1$. Hence, $w_2w_1(L^\infty(Y_1)\rtimes\Sigma_1)w_1^*w_2^*\subset L^\infty(Y_1)\rtimes\Sigma_1$. Next, note that the quasi-normalizer of $L^\infty(Y_1)\rtimes\Sigma_1$ inside $M$ equals itself by Lemma \ref{mixing}(1). Hence, $w_2w_1\in L^\infty(Y_1)\rtimes\Sigma_1$, which shows that $M_0:=L^\infty(Y_1)\rtimes\Sigma_1= w_2^*(L^\infty(Y_2)\rtimes\Sigma_2)w_2$.

Next, we note that $L^\infty(Y_1)$ is unitarily conjugate to $w_1^*L^\infty(Y_2)w_1$ inside $M_0.$ 
Indeed, by applying Lemma \ref{small} we deduce that $L^\infty(Y_1)\prec_{M_0} w_1^*L^\infty(Y_2)w_1$ and the claim follows by \cite[Theorem A.1]{Po01}.
Hence, there exists a $*$-isomorphism $\Psi$ of $M$ that sends 
$L^\infty(Y_1)\rtimes\Sigma_1$ onto $L^\infty(Y_2)\rtimes\Sigma_2$ and $L^\infty(Y_1)$ onto $L^\infty(Y_2)$. 

Suppressing $\Psi$ from the notation, we can assume that $M_0=L^\infty(Y_1)\rtimes\Sigma_1=L^\infty(Y_2)\rtimes\Sigma_2$ such that $L^\infty(Y_1)=L^\infty(Y_2)$. 
We continue with the following claim.

{\bf Claim 2.} There exist two unitaries $u_1,u_2\in M_0$ such that
$u_1(L^\infty(X_1)\rtimes\Sigma_1) u_1^*\subset L^\infty(X_2)\rtimes\Sigma_2$ and $u_2(L^\infty(X_2)\rtimes\Sigma_2)u_2^*\subset  L^\infty(X_1)\rtimes\Sigma_1.$

{\it Proof of Claim 2.} We will show only the first part. Let $\alpha_t$ be the free product deformation as in Definition \ref{def:free} associated to the coinduced action $\Gamma_2\car Y_2$.

Assume first that $\Sigma_1$ has property (T) or it is measure equivalent to a product of two non-amenable groups. As before, we derive that $\alpha_t\to id$ uniformly on $(L(\Sigma_1))_1$ by Lemma \ref{L:ME}. If $\Sigma_1\car (X_1,\mu_1)$ is rigid, then 
$\alpha_t\to id$ uniformly on $(L^\infty(X_1)\rtimes\Sigma_1)_1$. By applying the moreover part of Theorem \ref{rigid2} in the case $F_0=\{\Sigma_1\}$, we derive that there exists a unitary $u_1\in M_0$ such that
$u_1(L^\infty(X_1)\rtimes\Sigma_1) u_1^*\subset L^\infty(X_2)\rtimes\Sigma_2$ or $L(\Sigma_1)\prec_{M_0} L^\infty(Y_2)$. Note that the second possibility can not hold, since $\Sigma_1$ is infinite. If $\Sigma_1\car (X_1,\mu_1)$ is compact, we apply again \cite[Theorem 6.9]{Io08a} and Theorem \ref{rigid2} in the case $F_0=\{\Sigma_1\}$ and obtain that the claim holds.

Assume now that $\Sigma_1\car X_1$ is isomorphic to $\Sigma_1^1\times\Sigma_1^2\car X_1^1\times X_1^2$, with $\Sigma^1_1$ non-amenable. 
\cite[Lemma 2.6(2)]{DHI16} combined with a maximality argument imply that we can find a non-zero projection $z\in (L^\infty(X_1)\rtimes\Sigma_1)'\cap M$ such that 
$z((L^\infty (X_1^2)\rtimes\Sigma_1^2 )'\cap M)z$ is strongly non-amenable in $M$. Hence,
we can apply Theorem \ref{spectral gap} and deduce that
$\alpha_t\to id$ uniformly on $((L^\infty (X_1^2)\rtimes\Sigma_1^2 )z)_1$.  Lemma \ref{smaller} shows that $z\in (L^\infty(X_1)\rtimes\Sigma_1)'\cap M_0$. By applying the moreover part of Theorem \ref{rigid2} in the case $F_0=\{\Sigma_1\}$, we obtain that $L^\infty(X_1^2)\rtimes\Sigma_1^2\prec_{M_0} L^\infty(X_2)$ or $L^\infty(X_1)\rtimes\Sigma_1\prec_{M_0} L^\infty(X_2)\rtimes\Sigma_2$ or $L^\infty(X_1)\rtimes\Sigma_1\prec_{M_0}L^\infty(Y_2)$. Since $\Sigma_1^2$ is an infinite group, it follows that $L^\infty(X_1)\rtimes\Sigma_1\prec_{M_0} L^\infty(X_2)\rtimes\Sigma_2$. In this case, the last part of Theorem \ref{rigid2} implies actually the conclusion of the claim.
\hfill$\square$

Let $I_j=\Gamma_j/\Sigma_j$ and $Y^0_j=X_j^{I_j\setminus\{\Sigma_j\}}$, for both $1\leq j\leq 2$. Using these notations, we obtain that $L^\infty(X_1\times Y^0_1)\rtimes\Sigma_1= L^\infty(X_2\times Y^0_2)\rtimes\Sigma_2$. Note that by combining Lemma \ref{wm1} and Lemma \ref{wm2} it follows that the action $\Sigma_1\car Y_1^0$ is weakly mixing. Hence, we obtain by Lemma \ref{L:qn} that the quasi-normalizer of $L^\infty(X_1)\rtimes\Sigma_1$ inside $M_0$ equals itself. Next, we note that Claim 2 gives that $u_2u_1 (L^\infty(X_1)\rtimes\Sigma_1)u_1^*u_2^*\subset L^\infty(X_1)\times\Sigma_1$. Thus, $u_2u_1\in L^\infty(X_1)\rtimes\Sigma_1$ and therefore, $u_1(L^\infty(X_1)\rtimes\Sigma_1) u_1^*= L^\infty(X_2)\rtimes\Sigma_2$. In particular, it follows that $u_1L^\infty(X_1)u_1^*\prec_{M_0}L^\infty(X_2)$. Since $L^\infty(X_2)\subset M_0$ is regular, we obtain by Lemma \ref{smaller} that $u_1L^\infty(X_1)u_1^*\prec_{L^\infty(X_2)\rtimes\Sigma_2}L^\infty(X_2)$. Finally, we use \cite[Theorem A.1]{Po01} and obtain that $u_1L^\infty(X_1)u_1^*$ and $L^\infty(X_2)$ are unitarily conjugate in $L^\infty(X_2)\rtimes\Sigma_2$. By \cite{Si55,FM75} we get that $\Sigma_1\car X_1$ and $\Sigma_2\car X_2$ are orbit equivalent. 
\hfill$\blacksquare$

\end{document}